# ONE-DIMENSIONAL STEPPING STONE MODELS, SARDINE GENETICS AND BROWNIAN LOCAL TIME

By Richard Durrett and Mateo Restrepo

*Cornell University*

Consider a one-dimensional stepping stone model with colonies of size $M$ and per-generation migration probability $\nu$, or a voter model on $\mathbb{Z}$ in which interactions occur over a distance of order $K$. Sample one individual at the origin and one at $L$. We show that if $M\nu/L$ and $L/K^2$ converge to positive finite limits, then the genealogy of the sample converges to a pair of Brownian motions that coalesce after the local time of their difference exceeds an independent exponentially distributed random variable. The computation of the distribution of the coalescence time leads to a one-dimensional parabolic differential equation with an interesting boundary condition at 0.

**1. Introduction.** Cox and Durrett [6] and Zähle, Cox and Durrett [15] have recently studied the two-dimensional stepping stone model. Space is represented as a torus $\Lambda(L) = (\mathbb{Z} \bmod L)^2$. To avoid a factor of 2 and to make the dynamics easier to describe, we suppose that at each point $x \in \Lambda(L)$ there is a colony of $M$ haploid individuals labeled $1, 2, \ldots, M$. Each individual in the system is replaced at rate 1. With probability $1 - \nu$ it is replaced by a copy of an individual chosen from the same colony. If the individual is in colony $x$, then with probability $\nu$ it is replaced by a copy of one chosen from nearby colony $y \neq x$ with probability $q(y - x)$ where the difference is computed componentwise modulo $L$, and the representative of the equivalence class chosen from $(-L/2, L/2]^2$. Here $q(z)$ is an irreducible probability on $\mathbb{Z}^2$ with $q(0,0) = 0$, finite range and the same symmetry as $\mathbb{Z}^2$: $q(x_1, x_2) = q(-x_1, -x_2)$ and $q(x_1, x_2) = q(x_2, x_1)$. These assumptions imply that jumps according to $q$ have mean 0 and covariance $\sigma^2 I$.

When $M = 1$ the stepping stone model reduces to the voter model, but being able to consider colony size $M > 1$ enriches the behavior of the model.

As in the voter model, we can define a genealogical process for each individual that traces the source of its genetic material backward in time. For one individual this is a random walk that moves to a randomly chosen individual in the same colony with probability $1 - \nu$ and otherwise jumps to a new colony chosen according to $q$. The genealogies of two individuals are random walks that coalesce with probability $1/M$ on each jump when they land in the same colony. We will call $q$ the *dispersal distribution* since it is the jump distribution for the genealogical process. If the migration rate times the colony size, $M\nu$, is large enough, then the population behaves as a homogeneously mixing unit. Let $t_0$ be the coalescing time of two lineages and let $\pi$ denote that the two individuals are chosen at random from the population. Cox and Durrett [6] have shown

THEOREM 1.   *If $L \to \infty$ and $(2\pi\sigma^2)M\nu/\log L \to \alpha \in (0, \infty]$, then*
$$P_\pi\left(2t_0 > \frac{1+\alpha}{\alpha}ML^2 t\right) \to e^{-t}.$$

In genetics terms, the system behaves as a homogeneously mixing population of "effective" size $ML^2(1+\alpha)/\alpha$. As $\alpha \to \infty$ this converges to the actual population size, indicating that the critical size of $M\nu$ for interesting behavior is $O(\log L)$. One finds more interesting behavior when individuals are sampled from a $L^\beta \times L^\beta$ square of colonies, but those results are not relevant here, so we refer the reader to Cox and Durrett [6] and Zähle, Cox and Durrett [15] for details.

Here, we will be interested in investigating similar questions for the one-dimensional stepping stone model. Although we live in a two-dimensional world, this case is relevant for applications. Many species, such as sea lions and abalone, live along a coastline that is essentially one-dimensional. For example, Bowen and Grant [5] have studied sardines at five different sites in the Indian and Pacific oceans. Wilkins and Wakeley's [14] analysis of this data using the one-dimensional stepping stone model was the inspiration for this study.

Although the most natural setting to pursue our results would be a one-dimensional interval or a ring of colonies, we will, for technical reasons, study the stepping stone model on $\mathbb{Z}$. The setup is the same as that of Cox and Durrett [6] described above. There are $M$ haploid individuals per colony and nearest-neighbor migration occurs with probability $\nu$. We sample one individual from the colony at 0, and another from the colony at $L$. If $M = 1$, then the two lineages will coalesce the first time they enter the same colony. Our first question is how large should $M\nu$ need to be for the system to have more interesting behavior? Since migration occurs with probability $\nu$, it takes time $O(L^2/\nu)$ for the difference in the locations of the two lineages to change by $O(L)$. In this time the difference will visit a given value between 0 and $L$ an average of $L/\nu$ times, so if we want the probability of coalescence to be positive but not certain, this should be $O(M)$.

THEOREM 2. *Consider a one-dimensional stepping stone model with $M$ haploid individuals per colony and nearest-neighbor migration with probability $\nu$. Sample one individual from the colony at 0, and another from the colony at $L$. If $L \to \infty$ and $M\nu/L \to \alpha \in (0,\infty)$, then $2t_0/(L^2/\nu)$ converges in distribution to $\ell_0^{-1}(\alpha\xi)$, where $\ell_t(0)$ is the local time at 0 for a standard Brownian motion starting at 1, and $\xi$ is independent with a mean 1 exponential distribution.*

Note that as $\alpha \to 0$ the limit becomes the hitting time of 0 and that as $\alpha \to \infty$ the limit $\to \infty$.

We are, of course, not the first to have considered this problem. Writing things in our notation, Maruyama [12] considered a ring of $L$ colonies with $M$ diploid individuals per colony. He did not formulate his result as a limit theorem, but by filling in a few details in the Appendix, we can use his computations to show that if $M\nu/L \to \alpha$,

$$(1) \qquad E_0(\exp(-\lambda t_0/(L^2/\nu))) \to (1 + 4\alpha\sqrt{\lambda})^{-1}.$$

It would be interesting to derive this formula using a generalization of Theorem 2 to the circle, and computations for the local time at 0 of a Brownian motion on the circle.

Wilkins and Wakeley [14] modeled space as $\{0, 1/L, 2/L, \ldots, 1\}$ with one individual per site, and used a dispersal distribution that is a normal distribution with a small variance $\sigma^2$ with reflecting boundary conditions on the ends. They analyzed the system by simulation and numerical solution of differential equations for various combinations of $L$ and $\sigma^2$. Here we will consider the corresponding problem on $\mathbb{Z}$, sample one individual from 0 and one from $L$, and suppose dispersal distance is of order $K$. If the dispersal is nearest neighbor, the two lineages cannot cross each other without coalescing. To see how large $K$ has to be for the system to have interesting behavior, we note that it takes roughly $L^2/K^2$ jumps to move distance $L$, and at this point the difference between the two locations will have visited a typical value between 0 and $L$ about $L/K^2$ times. If we take $K = c\sqrt{L}$, then the expected number of visits to 0 converges to a positive finite limit, and the probability of coalescence is positive but not certain.

To state the result and to write its proof, it is convenient to introduce another parameter $N$ and let $K = N^{1/2}$ and $L = O(N)$. We make the following assumptions about the dispersal distribution $q^N$:

1. symmetry: $q^N(z) = q^N(-z)$,

2. the variance $\sum_{z \in \mathbb{Z}} z^2 q^N(z) = \sigma_N^2 N$ with $\sigma_N \to \sigma \in (0, \infty)$,
3. there is an $h > 0$, independent of $N$, so that $q^N(z) \geq h/\sqrt{N}$ for $|z| \leq N^{1/2}$,
4. exponential tails: $q^N(z) \leq C \exp(-c|z|/\sqrt{N})$.

These assumptions contain uniform, bilateral exponential and normal distributions as special cases. The last condition is strong but is convenient since it allows us to choose $B$ so that

$$\sum_{|z| \geq B\sqrt{N} \log N} q^N(z) \leq N^{-2}.$$

Since the limit theorem involves times of order $N$, we can suppose without loss of generality that

5. $q^N(z) = 0$ for $|z| > B\sqrt{N} \log N$,

since the probability of having a jump larger than $B\sqrt{N} \log N$ by time $N$ is $\leq 1/N$. The constant $B$ is special and the letter $B$ is reserved for its value. Here and in what follows, $c$ and $C$ are positive finite constants whose values are unimportant and will change from line to line, while $O(f(N))$ indicates a quantity that can be bounded by $Cf(N)$, with $C$ independent of $N$.

THEOREM 3. *Consider a sequence of voter models on $\mathbb{Z}$ with jumps at rate 1, and dispersal distributions $q^N$, satisfying assumptions 1–5. If the positive numbers $L_N$ have $L_N/(\sigma N) \to x_0 \geq 0$, then $2t_0/N$ converges in distribution to $\ell_0^{-1}(\sigma \xi/2)$, where $\ell_0$ is the local time at 0 of a standard Brownian motion started from $x_0$ and $\xi$ is independent with a mean 1 exponential distribution.*

Again, as $\sigma \to 0$ the limit becomes the hitting time of 0, and as $\sigma \to \infty$ the limit $\to \infty$.

To get a more explicit description of the distribution of the limits in Theorems 2 and 3 we would like to compute

$$P_x(\ell_0^{-1}(\xi/\lambda) > t) = P_x(\lambda \ell_0(t) < \xi) = E_x \exp(-\lambda \ell_0(t)).$$

Formula 1.3.7 in Borodin and Salaminen's [4] *Handbook of Brownian Motion* tells us that

$$E_x(e^{-\lambda \ell_0(t)}; W_t \in dz)$$

(2)
$$= \frac{1}{\sqrt{2\pi t}} e^{-(z-x)^2/2t} dz$$

$$- \frac{\lambda}{2} \exp((|z| + |x|)\lambda + \lambda^2 t/2) \operatorname{Erfc}\left(\frac{\lambda^2 \sqrt{t}}{\sqrt{2}} + \frac{|z| + |x|}{\sqrt{2t}}\right) dz,$$

where Erfc is the error function, that is, the upper tail of the normal distribution.

Another approach to computing $u(t,x) = E_x \exp(-\lambda \ell_0(t))$ is to note that for $x \neq 0$ it satisfies the heat equation

$$\frac{\partial u}{\partial t} = \frac{1}{2}\frac{\partial^2 u}{\partial x^2}.$$

To determine the boundary condition at 0, we run Brownian motion until $\tau_h = \inf\{t : B_t \notin (-h,h)\}$ and use symmetry $u(t,x) = u(t,-x)$ to conclude that

$$u(t,0) = E_0(e^{-\lambda \ell_0(\tau_h)} u(t - \tau_h, h); \tau_h \leq t) + P_0(\tau_h > t).$$

The strong Markov property implies that $\ell_0(\tau_h)$ is exponentially distributed. Let $D_\varepsilon(\tau_h)$ be the number of downcrossings of $(0,\varepsilon)$ by reflecting Brownian motion before it hits $h$. $D_\varepsilon(\tau_h)$ is geometrically distributed with mean $h/\varepsilon$ and $\lim_{\varepsilon \to 0} \varepsilon D_\varepsilon(t) = \ell_0(t)$ (see, e.g., page 48 of Itô and McKean [10]), so $E_0 \ell_0(\tau_h) = h$ and

$$E_0(e^{-\lambda \ell_0(\tau_h)}) = \frac{1/h}{\lambda + 1/h} = \frac{1}{1 + \lambda h}.$$

Using the explicit formula in (2) or the fact that $u(t,x)$ satisfies the heat equation with a bounded boundary condition on $[0,\infty) \times \{0\}$ shows $u(t,x)$ is Lipschitz continuous on $[0,T] \times [-K,K]$. Since $\tau_h$ has the same distribution as $h^2 \tau_1$, $|u(t - \tau_h, h) - u(t,h)| = O(h^2)$. Using this with $P_0(\tau_h > t) = o(h)$, we have

$$\begin{aligned}
\frac{\partial u}{\partial x}(t, 0+) &= \lim_{h \to 0} \frac{u(t,h) - u(t,0)}{h} \\
&= u(t,0) \lim_{h \to 0} \frac{1 - E_0(e^{-\lambda \ell_0(\tau_h)})}{h} = \lambda u(t,0).
\end{aligned}$$

The remainder of the paper is devoted to proofs. Theorem 2 is fairly straightforward to prove. Let $Z_t^N$ be the difference between the colony numbers for the two lineages, and let $Y_m^N$ be the embedded jump chain, which jumps when a lineage changes colonies. $Y^N(L^2 \cdot)/L$ converges to a Brownian motion starting from 1. Using the fact that $|B_t^0| - \ell_0(t)$ is a martingale, it is easy to show that if $V_m^N$ is the number of visits to 0 by $Y_m^N$ then $V^N(L^2 \cdot)/L$ converges to the local time $\ell_0$. (Borodin [3] proved this for aperiodic mean 0, finite-variance random walks.) Each visit to 0 by $Y_m^N$ brings a probability of coalescence of roughly $\nu/(\nu + 1/N)$ for our two lineages, and the result follows from routine calculations. See Section 2 for details.

It is easy to give an intuitive proof of Theorem 3 along similar lines. The difference in the location between two lineages in the genealogy of voter model is a continuous-time random walk that jumps at rate 2, so it is enough to consider the embedded discrete-time jump chain. Let $X_k^N$ be a random walk with jump distribution $q^N$. Let $1/2 < a < 1$. The number of visits to $I = [-N^a, N^a]$ by time $t$, divided by $2N^a$, converges to the local time at 0 of a Brownian motion. If we look at the chain $X_k^N$ only when it is in $I$, then we get a Markov chain that mixes more rapidly than its expected time to hit 0, so a result of Aldous and Fill [1] implies that the hitting time of 0 for the chain viewed on $I$ has approximately an exponential distribution.

To complete the proof outlined in the previous paragraph, one must prove that the excursions off of $I$ are sufficiently independent of the behavior in $I$ so that the exponential waiting time and the local time are asymptotically independent. We have not been able to formalize this intuition, so we will instead pursue an approach based on the downcrossing definition of local time. Let $T_0 = \inf\{k : |X_k^N| < N^{5/6}\}$ and for $m \geq 0$ let

$$S_m = \inf\{k > T_m : |X_k^N| > 2N^{5/6}\},$$
$$T_{m+1} = \inf\{k > S_m : |X_k^N| < |X_{S_m}^N| - N^{5/6}\}.$$

Visits to 0 can only occur during $[T_m, S_m]$, while most of the time is in the intervals $[S_m, T_{m+1}]$. The definition of $T_{m+1}$ is chosen so that the distribution of $T_{m+1} - S_m$ is independent of $X^N(S_m)$, and this allows us to get the desired asymptotic independence.

Section 3 gives the proof of Theorem 3 modulo three propositions that are established later. Let $M^N(n) = \sup\{m : S_m \leq n\}$ be the number of cycles completed by time $n$. Proposition 1, proved in Section 4, gives the convergence of $M^N(Nt)/N^{1/6}$ to local time. Proposition 2, proved in Section 5, shows that the time spent in the intervals $[T_m, S_m]$ is a small fraction of the total time. Proposition 3 gives asymptotics for the probability of hitting 0 before time $S_m$ for the possible values of $X^N(T_m)$, which are $\pm N^{5/6} + O(N^{1/2} \log N)$. Proposition 3 is the most difficult part of the proof. It relies on estimates for the potential kernel, which are based on results for the Green's function, which in turn come from a local central limit theorem. The technical problem is that all of our estimates must be uniform in $N$. These details occupy Sections 6 and 7.

**2. Proof of Theorem 2.** Let $Z_t^N$ be the difference in the colony numbers at time $t$. Let $Y_m^N$ be the discrete-time embedded chain that jumps whenever one of the two lineages changes colonies, and continues jumping even after the two lineages have coalesced. $Y_m^N$ is a simple random walk. Recalling $Y_0^N = L$, we let $W^N(t) = Y_{[L^2 t]}^N / L$. Since $W^N$ converges in distribution to a standard Brownian motion $W(\cdot)$, $\mathcal{C} = C([0, \infty), \mathbb{R})$ with the topology of uniform convergence on compact time intervals is a complete separable metric space, Skorokhod's theorem implies that we can assume these processes have been constructed on the same space so that $W^N(\cdot) \to W(\cdot)$ almost surely. See, for example, Theorem 3.3 on page 7 of Billingsley [2].

Let $V_m^N$ be the number of visits to 0 by $Y_k^N$, $k \leq m$. The next result has been proved for finite-variance random walks by Borodin [3]. To keep this paper self-contained, we will give a simple proof for the nearest-neighbor case.

LEMMA 1. $V^N(L^2 \cdot)/L \to \ell_0(\cdot)$, the local time at 0 for $W$, almost surely in $\mathcal{C}$.

PROOF. Let $A^N(t) = V^N(L^2 \cdot)/L$. An easy computation for simple random walk shows that for any stopping time $S$

$$E_x|A^N(S+t) - A^N(S)| \leq E_0|A^N(t)| \leq \frac{1}{L}\left(1 + \sum_{k=1}^{L^2 t} C/\sqrt{k}\right) \leq C\sqrt{t},$$

so by Aldous' criterion (see, e.g., Theorem 4.5 on page 320 of Jacod and Shiryaev [10]) the sequence $A^N$ is tight. Let $A^{N_k}$ be a convergent subsequence with limit $A$. $|W^{N_k}(t)| - A^{N_k}(t)$ is a martingale. Using the $L^2$ maximal inequality on the random walk, and the dominated convergence theorem on the increasing process, both processes converge to their limits in $L^1$. Since conditional expectation is a contraction in $L^1$, it follows that $|W(t)| - A(t)$ is a martingale. $\ell_0(t)$ is the increasing process associated with $|W(t)|$. See, for example, (11.2) on page 84 of Durrett [7]. By the uniqueness of the Doob–Meyer decomposition $A(t) = \ell_0(t)$. This shows that there is only one subsequential limit, so the entire sequence converges to $\ell_0(t)$. □

To move this result from $Y^N$ to $Z^N$, we note that time $m$ in $Y^N$ corresponds to a time $\sim m/2\nu$ in $Z^N$, and hence time $L^2 t/2\nu$ in $Z^N$ corresponds to a time $\sim L^2 t$ in $Y^N$, where as usual $a_N \sim b_N$ means $a_N/b_N \to 1$. Now $Z^N$ will have a geometric number of chances with mean $1/\nu$ for coalescence between jumps of $Y^N$ so the probability of no coalescence is

$$\sum_{j=1}^{\infty}(1-\nu)^{j-1}\nu(1-1/N)^j = \frac{\nu(1-1/N)}{1/N + \nu(1-1/N)} \sim \frac{N\nu}{1+N\nu}.$$

Recall our assumptions imply $N\nu \to \infty$ and hence $N \to \infty$.

When $m = L^2 t$, the number of visits to 0 by $Y_m^N$ will be $\sim L\ell_0(t)$ and hence the probability of no coalescence is

$$= \left(1 - \frac{1}{N\nu}\right)^{L\ell_0(t)} \to e^{-(1/\alpha)\ell_0(t)}.$$

If $\xi$ is a mean 1 exponential, the right-hand side can be written as

$$P((1/\alpha)\ell_0(t) < \xi) = P(\ell_0^{-1}(\alpha\xi) > t),$$

which completes the proof of Theorem 2.

**3. Proof of Theorem 3.** Here we give the proof, assuming the truth of three propositions that will be proved in the next three sections. Let $X_k^N$, $k \geq 0$, be a discrete-time random walk with jump distribution $q^N$. To avoid some annoying little details, it is convenient to suppose that $X_0^N = x_N \geq 2N^{5/6}$. To extend to the general case, it is enough to show that starting from $x \neq 0$ the probability of hitting 0 before time $S_0$ defined below tends to 0, but this follows from Lemma 5.

Define two interleaved sequences of stopping times as follows. Let $T_0 = -1$ and for $m \geq 0$ let

$$S_m = \inf\{k > T_m : |X_k^N| > 2N^{5/6}\},$$
$$T_{m+1} = \inf\{k > S_m : |X_k^N| < |X_{S_m}^N| - N^{5/6}\}.$$

$S_m$ is the exit time from the larger strip $[-2N^{5/6}, 2N^{5/6}]$. Since

$$2N^{5/6} - BN^{1/2}\log N \leq |X^N(S_m)| \leq 2N^{5/6} + BN^{1/2}\log N,$$

$T_m$ is almost the hitting time of the smaller strip $[-N^{5/6}, N^{5/6}]$. The advantage of this definition is that the processes $\{|X^N(S_m + k)| - |X^N(S_m)|, 0 \leq k \leq T_m - S_m\}$ are identically distributed for $m \geq 0$ and independent of $\mathcal{F}(S_m)$. Here and in what follows, we will write $X^N(S_m)$ instead of $X_{S_m}^N$ to avoid double subscripts.

Let $M^N(n) = \sup\{m : S_m \leq n\}$ be the number of cycles completed by time $t$ and let $L^N(n) = |\{1 \leq m \leq M^N(n) : X^N(S_{m-1})X^N(S_m) < 0\}|$ be the number of crossings of $[-2N^{1/6}, 2N^{1/6}]$ by the random walk. Our first result to be proved later is:

PROPOSITION 1. *Suppose $x_N/\sigma N \to x_0$. Then*

$$2L^N(Nt)/N^{1/6} \Rightarrow \sigma\ell_0(t) \quad and \quad M^N(Nt)/2N^{1/6} \Rightarrow \sigma\ell_0(t),$$

*where $\ell_0(\cdot)$ is the local time at 0 of a standard Brownian motion started at $x_0$.*

Let $J = \inf\{m : \exists k \in [T_m, S_m], X_k^N = 0\}$. The fact that one-dimensional finite-range random walks are recurrent implies $J < \infty$. By the definitions of $S_m$ and $T_m$, $t_0 \in [T_J, S_J]$. Splitting things up according to the value of $J$,

$$\sum_{j=0}^{\infty} P\{J = j, T_j > Nt\} \leq P\{t_0 > Nt\} \leq \sum_{j=0}^{\infty} P\{J = j, S_j > Nt\}.$$

We will show that both series converge to the same limit as $N \to \infty$, thereby proving that $P\{t_0 > Nt\}$ converges to this limit as well. We first truncate the sums by neglecting the terms having $j > N^{2/9}$:

$$(3) \quad \sum_{j=0}^{N^{2/9}} P\{J = j, T_j > Nt\} \leq P\{t_0 > Nt\} \leq \sum_{j=0}^{N^{2/9}} P\{J = j, S_j > Nt\} + \varepsilon_1^N,$$

where, as Lemma 2 will show,

$$\varepsilon_1^N \leq \sum_{j=N^{2/9}}^{\infty} P\{J=j\} = O(e^{-cN^{1/18}}).$$

Defining $A_j = T_0 + \sum_{m=1}^{j}(T_m - S_{m-1})$ and $B_j = \sum_{m=0}^{j}(S_m - T_m)$, we can write $S_j = A_j + B_j$. For the reader's intuition, we note that $T_m - S_{m-1}$ is the hitting time of a half-line, while $S_m - T_m$ is the exit time from a bounded strip. The first variable has infinite mean and the latter finite variance, so we expect $A_j \gg B_j$ for large $j$.

It is clear that

(4) $$\sum_{j=0}^{N^{2/9}} P\{J=j, A_j > Nt\} \leq \sum_{j=0}^{N^{2/9}} P\{J=j, T_j > Nt\}.$$

Our next task is to argue that

(5) $$\sum_{j=0}^{N^{2/9}} P\{J=j, S_j > Nt\} \leq \sum_{j=0}^{N^{2/9}} P\{J=j, A_j > Nt - 2N^{17/18}\} + \varepsilon_2^N,$$

where $\varepsilon_2^N$ is another small error, this time $O(N^{-1/3})$. To prove the last inequality we note that, for any $j \leq N^{2/9}$,

$$\{J=j, S_j > Nt\} \subset \{J=j, A_j > Nt - 2N^{17/18}\} \cup \{J=j, B_{N^{2/9}} > 2N^{17/18}\}.$$

In the last equality we should have written the integer part $[N^{2/9}]$, but in what follows we will ignore these insignificant details. Taking now the union over $j$, (5) will follow from the following proposition.

PROPOSITION 2. *For large $N$, $P\{B_{N^{2/9}} > 2N^{17/18}\} \leq CN^{-1/3}$.*

Combining inequalities (3), (4) and (5), we can restrict ourselves to estimating probabilities of the form $P\{J=j, A_j > Ns\}$. The two events here are almost independent. $A_j$ is determined by the behavior of increments of the random walk in the intervals $[S_m, T_{m+1}]$, while $J$ is determined by the behavior in $[T_m, S_m]$. There is some dependence that comes through the value of the starting points $X^N(T_m)$, but because of assumption 5, these are all within distance $BN^{1/2} \log N$ of $N^{5/6}$ or $-N^{5/6}$. As the reader can probably guess, the variability in the starting point makes little difference:

PROPOSITION 3. *Suppose $|x - N^{5/6}| \leq BN^{1/2} \log N$ and let $H_I^N(x,0)$ denote the probability that the random walk $X^N$ started at $x$ hits $0$ before leaving the set $I = [-2N^{5/6}, 2N^{5/6}]$. There is a constant $C$ so that*

$$\left| N^{1/6} H_I^N(x,0) - \frac{1}{\sigma^2} \right| \leq CN^{-1/6} \log N.$$

The bound in Proposition 3 is uniform over the possible values of $X^N(T_m)$, so for simplicity we will write $c_N = 1/\sigma^2 + O(N^{-1/6}\log N)$.

LEMMA 2. *For every $u > 0$ and $j$ we have*

(6) $$P\{J=j, A_j > u\} = \frac{c_N}{N^{1/6}}\left(1 - \frac{c_N}{N^{1/6}}\right)^j P\{A_j > u\},$$

*and hence* $\sum_{j=N^{2/9}}^{\infty} P\{J=j\} \leq (1 - c_N/N^{1/6})^{N^{2/9}} \leq \exp(-cN^{1/18})$.

PROOF. Let $I_m = 1\{t_0 \in [T_m, S_m]\}$ and let $\Delta_k = A_k - A_{k-1}$ for $k \geq 0$, where $A_{-1} = 0$. Using the strong Markov property, Proposition 3, the fact that $\Delta_j$ is independent of $\mathcal{F}(S_{j-1})$ and induction, it is easy to see that

$$P(\Delta_0 = v_0, I_0 = 0, \Delta_1 = v_1, \ldots, I_{j-1} = 0, \Delta_j = v_j, I_j = 1)$$
$$= \frac{c_N}{N^{1/6}}\left(1 - \frac{c_N}{N^{1/6}}\right)^j \prod_{k=0}^{j} P(\Delta_k = v_k).$$

Since the $\Delta_k$ are independent, the desired result follows by summing over $v_0, \ldots, v_k$ that sum to more than $u$. □

The lower bound in (4) and the upper bound in (5) are similar, so it is enough to investigate the lower bound. Using Lemma 2 on the left-hand side of (4) gives

$$P\{t_0 > Nt\} \geq \sum_{j=0}^{N^{2/9}} \frac{c_N}{N^{1/6}}\left(1 - \frac{c_N}{N^{1/6}}\right)^j P\{A_j > Nt\}.$$

Using Proposition 2, we get

$$P\{t_0 > Nt\} \geq \sum_{j=0}^{N^{2/9}} \frac{c_N}{N^{1/6}}\left(1 - \frac{c_N}{N^{1/6}}\right)^j P\{S_j > Nt + 2N^{17/18}\} + \varepsilon_3^N,$$

where $\varepsilon_3^N$ is an error of order $N^{-1/3}$. Recalling the definition of $M^N$, the above is

(7) $$= \sum_{j=0}^{N^{2/9}} \frac{c_N}{N^{1/6}}\left(1 - \frac{c_N}{N^{1/6}}\right)^j P\{M^N(Nt + 2N^{17/18}) < j\}.$$

Proposition 1 implies that
$$P(M^N(Nt + 2N^{17/18}) < sN^{1/6}) \to P(\sigma \ell_0(t) < s/2).$$

Let $c_0 = 1/\sigma^2 = \lim_{N\to\infty} c_N$. The dominated convergence theorem now implies that (7) converges to

$$\int_0^{\infty} c_0 e^{-c_0 s} P\{\sigma \ell_0(t) < s/2\}\, ds.$$

Introducing a mean 1 exponential random variable, $\xi$, independent of $L_0(t)$, and recalling $c = 1/\sigma^2$ is 1 over the mean of the exponential, this can be written as

$$P\{\sigma \ell_0(t) < \sigma^2 \xi/2\} = P\{\ell_0^{-1}(\sigma \xi/2) > t\},$$

which is the conclusion of Theorem 3. It remains to prove the three propositions.

**4. Proof of Proposition 1.** Consider the sequence of random walks $X_k^N = x_N + \sum_{i=1}^k \xi_i^N$ where for each $N$, $x_N \geq 2N^{5/6}$ and the variables $\xi_i^N$, $i \geq 1$, are i.i.d. with distribution $q^N$. Define the sequence of stopping times $K_0 = 0$ and for $j \geq 0$

(8)
$$K_{2j+1} = \inf\{k > K_{2j} : X_k^N < -2N^{5/6}\},$$
$$K_{2j+2} = \inf\{k > K_{2j+1} : X_k^N > 2N^{5/6}\}.$$

In words, the $K_{2j+1}$ correspond to times at which the random walk finishes a down crossing of the interval $[-2N^{5/6}, 2N^{5/6}]$ and the $K_{2j+2}$ correspond to times at which the random walk finishes an up crossing of the same interval.

To connect with the definitions given just before Proposition 1 in the previous section, note that $\{S_m : m \geq 0\} \supset \{K_k : k \geq 0\}$ (it is for this reason that we want $x_N \geq 2N^{5/6}$), so we have

$$L^N(n) = \sup\{j : K_j \leq n\}.$$

Here and in what follows, even though $\sigma_N \to \sigma$ we will drop the subscript $N$ for simplicity.

LEMMA 3. *Suppose that $x_N/\sigma N \to x_0$, the $\xi_i^N$ are i.i.d. with $\mathrm{E}\xi_i^N = 0$, $\mathrm{E}(\xi_i^N)^2 = N\sigma^2$ and $\mathrm{E}(\xi_i^N)^4 \leq CN^2$. Then*

$$2N^{-1/6} L^N([Nt]) \Rightarrow \sigma \ell_0(t)$$

*as $N \to \infty$, where $\ell_0(t)$ denotes the local time at 0 for a standard Brownian motion starting from $x_0$.*

PROOF. We first rescale the random walks by letting

$$S_k^N = \frac{X_k^N}{N\sigma} = \frac{1}{\sigma\sqrt{N}} \sum_{i=1}^k \frac{\xi_i^N}{\sqrt{N}}.$$

Let $Y^N(t) = S_{[Nt]}^N$. Our first task is to argue that it is possible to define the $Y^N$'s and a Brownian motion $B$ on the same probability space $\Omega$, so that for each fixed $t$, the events

(9) $$\Omega_N = \left\{ \sup_{0 \leq s \leq t} |B(s) - Y^N(s)| \leq N^{-5/24} \right\}$$

satisfy $P(\Omega_N) \to 1$.

To prove this, we begin by recalling a well-known construction of Skorohod, see, for example, Section 7.6 in Durrett [7]. Given a Brownian motion $B$ and a value of $N$, this procedure constructs a sequence of stopping times $T_k^N$, $k \geq 1$, that satisfy

$$B(T_k^N) \stackrel{d}{=} S_k^N = Y^N(k/N)$$

and are such that the increments $\tau_i^N = T_i^N - T_{i-1}^N$ are independent, non-negative random variables having mean $E\tau_i^N = E(\xi_i^N/\sigma N)^2 = 1/N$, and variance

$$\text{var}(\tau_i^n) \leq CE(\xi_i^N/\sigma N)^4 \leq C/N^2.$$

For $s \in [k/N, (k+1)/N)$, we have

$$|Y^N(s) - B(s)| = |Y^n(k/N) - B(s)| \leq |B(T_k^N) - B(k/N)| + |B(k/N) - B(s)|.$$

We now fix $t$, and argue that there are sets $\Omega_N^1$ with $P(\Omega_N^1) \to 1$ on which

(10) $\qquad |T_k^N - k/N| < N^{-11/24} \qquad \text{for all } k \leq Nt.$

Kolmogorov's $L^2$ maximal inequality (see, e.g., (4.3) in Chapter 4 of Durrett [8]) applied to the martingale $T_k^N - k/N$ gives

$$P\left(\sup_{k \leq Nt}\left|T_k^N - \frac{k}{N}\right| \geq N^{-11/24}\right) \leq N^{11/12} Nt\, \text{var}(\tau_i^N) \leq CtN^{-1/12}.$$

By Lévy's result on the modulus of continuity for Brownian motion we can find sets $\Omega_N^2$, with $P(\Omega_N^2) \to 1$ and such that, on $\Omega_N^2$, $x, y \leq t$ and $|x - y| \leq N^{-11/24}$ imply (see, e.g., (4.10) in Chapter 7 of Durrett [8]),

$$|B(x) - B(y)| \leq 10(|x-y|\log(|x-y|^{-1}))^{1/2} \leq (1/2)|x-y|^{5/11},$$

the last inequality holding for large $N$ since $5/11 < 1/2$. On $\Omega_N^1 \cap \Omega_N^2$ we have for $s \in [k/N, (k+1)/N)$

$$|Y^N(s) - B(s)| \leq \left|B(T_k^N) - B\left(\frac{k}{N}\right)\right| + \left|B\left(\frac{k}{N}\right) - B(s)\right|$$
$$\leq \frac{1}{2}(N^{-(11/24)(5/11)} + N^{-5/11}) \leq N^{-5/24},$$

which proves (10).

Having established (10), the rest of the proof of Lemma 3 is straightforward. Let $a_N = (1/\sigma)2N^{-1/6}$ and $b_N = (1/\sigma)N^{-5/24}$. Using definitions similar to the $K_j$ in (8), we can define $\mathcal{L}_N^-(t)$ and $\mathcal{L}_N^+(t)$ to be the number of times the Brownian motion $B_s$, $0 \leq s \leq t$, has crossed the strips

$[-a_N + b_N, a_N - b_N]$ and $[-a_N - b_N, a_N + b_N]$, respectively. On the events $\Omega_N$ we have
$$\mathcal{L}_N^+(t) \leq L^N(Nt) \leq \mathcal{L}_N^-(t).$$
On the other hand, a classical result obtained by Lévy on the convergence of downcrossings to local time (see Itô and McKean [9], page 48) implies that, as $N \to \infty$,
$$(a_N + b_N)\mathcal{L}_N^+(t) \to \ell_0(t), \qquad (a_N - b_N)\mathcal{L}_N^-(t) \to \ell_0(t).$$
To check the constant, recall that one multiplies the number of downcrossings by the width of the strip, but here we count up- and downcrossings. This completes the proof of Lemma 3. □

To prove the convergence result for $M^N$ given in Proposition 1, we let
$$\Gamma(n) = |\{1 \leq m \leq n : X^N(S_{m-1})X^N(S_m) < 0\}|$$
and note that $L_n^N = \Gamma(M_n^N)$. Let $\gamma_m = 1$ if $X^N(S_{m-1})X^N(S_m) < 0$. We have
$$-B\sqrt{N}\log N \leq |X^N(T_m)| - N^{5/6} \leq B\sqrt{N}\log N,$$
so using the fact that $X_k^N$ is a martingale,
$$P(\gamma_m = 1 | \mathcal{F}(T_m)) = 1/4 + O(N^{-1/2}\log N). \tag{11}$$
Let $\bar{\Gamma}(n) = \Gamma(n) - \sum_{m=1}^n P(\gamma_m = 1|\mathcal{F}(T_m))$. $\bar{\Gamma}(n)$ is a martingale so the $L^2$ maximal inequality and the orthogonality of martingale increments imply
$$E\left(\sup_{m \leq n} \bar{\Gamma}(m)\right)^2 \leq C \sum_{m=1}^n E(\gamma_m - P(\gamma_m = 1|\mathcal{F}(T_m)))^2 \leq Cn.$$
Chebyshev's inequality implies
$$P\left(\sup_{m \leq n} \bar{\Gamma}(m) > n^{2/3}\right) \leq n^{-1/3}.$$
The last result when combined with (11) implies that with high probability
$$\Gamma(n) = n/4 + O(nN^{-1/2}\log N) + O(n^{2/3}).$$
We want to conclude from this that
$$L_n^N = \Gamma(M_n^N) \sim M_n^N/4.$$
To deal with the random index, we take $n = N^{1/5}$ and let $R = \inf\{r : M^N(r) \geq N^{1/5}\}$ to get
$$P\left(\sup_{s \leq Nt \wedge R} |L^N(s) - M^N(s)/4| > 2N^{2/15}\right) \leq N^{-1/15}.$$
Since $P(L^N(Nt) \geq N^{1/5}/5) \to 0$ by Lemma 3, we must have $P(R \leq Nt) \to 0$ and the proof of Proposition 1 is complete.

**5. Proof of Proposition 2.** In this section we will show that for large $N$,
$$P\{B_{N^{2/9}} > 2N^{17/18}\} \le CN^{-1/3},$$
where $B_j = \sum_{m=0}^{j} S_m - T_m$. To do this we will compute the mean and variance of $B_j$ and then use Chebyshev's inequality. For this we first need to compute the first two moments of $\eta_m = S_m - T_m$. If we assume $X^N(T_m) = N^{5/6}$, $|X^N(S_m)| = 2N^{5/6}$, and replace our random walk by a Brownian motion $B_t$ with variance $\sigma^2 Nt$, this would be easy. $B_t^2 - \sigma^2 Nt$ and $B_t^4 - 6\sigma^2 N B_t^2 t + 3\sigma^4 N^2 t^2$ are martingales so if $B_0 = N^{5/6}$ and $\eta = \inf\{B_t \notin [-2N^{5/6}, 2N^{5/6}]\}$, then using $|B_\eta| = 2N^{5/6}$ we have
$$4N^{10/6} - \sigma^2 N E\eta = N^{10/6},$$
$$16N^{20/6} - 6\sigma^2 N \cdot 4N^{10/6} E\eta + 3\sigma^4 N^2 E\eta^2 = N^{20/6}.$$
To prove this one must use the optional stopping theorem at $\eta \wedge m$ and then let $m \to \infty$. The details of using the monotone and dominated convergence theorem to justify the equalities are left to the reader. Solving gives
$$E\eta = 3N^{2/3}/\sigma^2,$$
$$E\eta^2 = 19N^{4/3}/\sigma^4.$$

These facts are approximately true for the random walk. We begin with the martingales. To compare with the previous calculation, recall that for the normal distribution $E\xi^4 = 3(E\xi^2)^2$.

LEMMA 4. *Suppose $X_k = X_0 + \xi_1 + \cdots + \xi_k$ where $E\xi_i = 0$, $E\xi_i^2 = \alpha$, $E\xi_i^3 = 0$ and $E\xi_i^4 = \beta$. Then $X_k^2 - k\alpha$ and*
$$X_k^4 - 6\alpha X_k^2 k + 3\alpha^2 k^2 + (3\alpha^2 - \beta)k$$
*are martingales.*

PROOF. The martingale $X_k^2 - k\alpha$ is well known. See, for example, Exercise 2.6 on page 235 of Durrett [8]. To check the second, expand $(X_k + \xi_{k+1})^4$ and use $E\xi_k = 0$ and $E\xi_k^3 = 0$ to conclude
$$E(X_{k+1}^4|\mathcal{F}_k) = X_k^4 + 6X_k^2 \alpha + \beta$$
and hence
$$E(X_{k+1}^4 - 6X_k^2(k+1)\alpha - \beta(k+1)|\mathcal{F}_k) = X_k^4 - 6X_k^2 k\alpha - \beta k.$$
To get the martingale we want, the $X_k^2$ on the left should be $X_{k+1}^2$. To correct this we note
$$E(-6(X_{k+1}^2 - X_k^2)(k+1)\alpha|\mathcal{F}_k) + 3\alpha^2(k+1)^2 + 3\alpha^2(k+1)$$
$$= -6\alpha^2(k+1) + 3\alpha^2(k+1)^2 + 3\alpha^2(k+1)$$
$$= 3\alpha^2[(k+1)^2 - (k+1)] = 3\alpha^2 k^2 + 3\alpha^2 k.$$

Adding the last two equations gives the desired result. □

In our case $\alpha = \sigma^2 N$ and $\beta \leq CN^2$. Letting $\mathcal{G}_{m-1} = \mathcal{F}(T_m)$ and using the optional stopping theorem on our first martingale with $|X^N(T_m)| \geq N^{5/6} - BN^{1/2}\log N$ and $|X^N(S_m)| \leq 2N^{5/6} + BN^{1/2}\log N$, we have

$$\sigma^2 N E(\eta_m | \mathcal{G}_{m-1}) \leq (2N^{5/6} + BN^{1/2}\log N)^2 - (N^{5/6} - BN^{1/2}\log N)^2$$
$$= 3N^{10/6} + O(N^{8/6}\log N),$$

and it follows that if $C_1 > (3/\sigma^2)$, then for large $N$

(12) $$E(\eta_m | \mathcal{G}_{m-1}) \leq C_1 N^{2/3}.$$

From the second martingale we get

$$E(X^N(S_m)^4 - 6\alpha X^N(S_m)^2 \eta_m + 3\alpha^2 \eta_m^2 + (3\alpha^2 - \beta)\eta_m | \mathcal{G}_{m-1}) = E(X^N(T_m)^4).$$

Rearranging and using $X^N(S_m)^4 \geq X^N(T_m)^4$ gives

$$3\alpha^2 E(\eta_m^2 | \mathcal{G}_{m-1}) \leq E(6\alpha X^N(S_m)^2 \eta_m - (3\alpha^2 - \beta)\eta_m | \mathcal{G}_{m-1}).$$

Using $|X^N(T_m)| \leq N^{5/6} + BN^{1/2}\log N$ and $|X^N(S_m)| \leq 2N^{5/6} + BN^{1/2}\log N$ with (12), $\alpha = \sigma^2 N$ and $\beta \leq CN^2$, gives

$$3\sigma^4 N^2 E(\eta_m^2 | \mathcal{G}_{m-1}) \leq [6(\sigma^2 N)(2N^{5/6} + BN^{1/2}\log N)^2 + CN^2] \cdot C_1 N^{2/3}.$$

The first term in the square brackets is of order $N \cdot N^{10/6} \gg N^2$. It follows that if $C_2 > 8C_1/\sigma^2$, then for large $N$

(13) $$E(\eta_m^2 | \mathcal{G}_{m-1}) \leq C_2 N^{4/3}.$$

To estimate the size of $B_j$, recall $\mathcal{G}_{m-1} = \mathcal{F}(T_m)$ for $m \geq 0$ and write

$$B_j = \sum_{m=0}^{j} E(\eta_m | \mathcal{G}_{m-1}) + \sum_{m=0}^{j} \eta_m - E(\eta_m | \mathcal{G}_{m-1}).$$

By (12), if $j \leq N^{2/9}$, then the first sum

$$\sum_1 \leq (j+1) C_1 N^{2/3} \leq 2C_1 N^{8/9}.$$

To bound the second sum, we use the orthogonality of martingale increments and (13) to conclude

$$E\left(\sum_2\right)^2 = \sum_{m=0}^{j} E(\eta_m - E(\eta_m | \mathcal{G}_{m-1}))^2 \leq (j+1) C_2 N^{4/3}.$$

When $j \leq N^{2/9}$, the right-hand side is $\leq 2C_2 N^{14/9}$:

$$P\left(\sum_2 \geq N^{17/18}\right) \leq 2C_2 N^{14/9} N^{-17/9} = C_2 N^{-1/3}.$$

Combining the bounds on $\sum_1$ and $\sum_2$ gives the conclusion of Proposition 2.

**6. Proof of Proposition 3.** Recall that the recurrent potential kernel is defined by

$$a(x,y) = \sum_{n=0}^{\infty}(p_n(x,y) - p_n(y,y)),$$

where $p_n$ is the $n$-step transition probability of the random walk. To see the reason for this definition, note that

(14)
$$\sum_x p(z,x)a(x,y) = \sum_{n=0}^{\infty}(p_{n+1}(z,y) - p_n(y,y)) = a(z,y), \qquad z \neq y,$$

$$\sum_x p(y,x)a(x,y) = \sum_{n=0}^{\infty}(p_{n+1}(y,y) - p_n(y,y)) = -1,$$

so $a$ is the analogue of the Green's function for recurrent random walks. The key to the proof of Proposition 3 is the following result whose proof is given in the next section. Let $\delta(x,y) = 1$ if $x = y$ and 0 otherwise.

PROPOSITION 4. *Assume a sequence of random walks satisfies assumptions 1–5 of Section 1. There is a constant $C$ independent of $N$ such that, for all $x$, their recurrent potential kernels satisfy*

$$\left|a^N(x,y) - \left(-1 + \delta(x,y) - \frac{|x-y|}{\sigma^2 N}\right)\right| \leq \frac{C}{\sqrt{N}}.$$

This estimate is only useful for $|x| \gg \sqrt{N}$. Our interest in this result is that it gives the following estimate on the Green's function $G_I^N(x,y)$, which is defined to be the expected number of visits to $y$ starting at $x$ before leaving the set $I$. If we let $\tau_I$ be the exit time from $I$, then in symbols,

$$G_I^N(x,y) = \sum_{k=0}^{\infty} P_x\{X_k^N = y, k < \tau_I\}.$$

We will be interested in the case $I = [-M, M]$ with $M = 2N^{5/6}$.

PROPOSITION 5. *There is a $C$ independent of $N$ such that for all $x$ and $y$*

$$\left|G_I^N(x,y) - \left(\delta(x,y) + \frac{M}{\sigma^2 N}\left[-\frac{|x-y|}{M} + \left(1 - \frac{xy}{M^2}\right)\right]\right)\right| \leq CN^{-1/3}\log N.$$

REMARK. To see that the formula in square brackets is reasonable, note that it vanishes when $x = M$ or $x = -M$ and for fixed $y$ is linear for $x \in [-M, y]$ and $x \in [y, M]$.

PROOF OF PROPOSITION 5. The first step is to note that (14) implies $a^N(X_n, y) + \sum_{m=0}^{n-1} \delta(X_m, y)$ is a martingale, so

$$G_I^N(x, y) = E_x[a^N(x, y) - a^N(X_{\tau_I}, y)]. \tag{15}$$

From (15) we have, for each fixed $N$, and $x, y \in [-M, M]$:

$$\begin{aligned}G_I^N(x, y) = a^N(x, y) &- P^x\{X_{\tau_I}^N > M\} E^x[a(X_{\tau_I}^N, y) | X_{\tau_I}^N > M] \\ &- P^x\{X_{\tau_I}^N < -M\} E^x[a(X_{\tau_I}^N, y) | X_{\tau_I}^N < -M]. \end{aligned} \tag{16}$$

Using now that $0 \leq X_{\tau_I}^N - M \leq BN^{1/2} \log N$ when $X_{\tau_I}^N > M$, the corresponding inequality for exiting at $-M$, and the fact that the random walk is a martingale, we have

$$x \leq P_x\{X_{\tau_I}^N > M\}(M + B\sqrt{N} \log N) + (1 - P_x\{X_{\tau_I}^N > M\})(-M),$$
$$x \geq P_x\{X_{\tau_I}^N > M\}M + (1 - P_x\{X_{\tau_I}^N > M\})(-M - B\sqrt{N} \log N).$$

Using these equations we have

$$\frac{M + x}{2M + B\sqrt{N} \log N} \leq P_x\{X_{\tau_I}^N > M\} \leq \frac{M + x + B\sqrt{N} \log N}{2M + B\sqrt{N} \log N},$$

and it follows that

$$P_x\{X_{\tau_I}^N > M\} = \frac{M + x}{2M} + O(\sqrt{N} \log N/M).$$

Subtracting from 1,

$$P_x\{X_{\tau_I}^N < -M\} = \frac{M - x}{2M} + O(\sqrt{N} \log N/M).$$

Using the last two formulas and Proposition 4 in (16),

$$\begin{aligned}G_I^N(x, y) = {}&-1 + \delta(x, y) - \frac{|x - y|}{\sigma^2 N} + O(1/\sqrt{N}) \\ &+ \left(\frac{M + x}{2M} + O(\sqrt{N} \log N/M)\right)\left(1 + \frac{M - y}{\sigma^2 N} + O(\log N/\sqrt{N})\right) \\ &+ \left(\frac{M - x}{2M} + O(\sqrt{N} \log N/M)\right)\left(1 + \frac{M + y}{\sigma^2 N} + O(\log N/\sqrt{N})\right). \end{aligned}$$

The worst error term is $O(\sqrt{N} \log N/M) = O(N^{-1/3} \log N)$. Ignoring the error terms, the sum of the second and third lines is

$$1 + \frac{M}{\sigma^2 N} - \frac{yx}{M\sigma^2 N}.$$

Adding this to the first line completes the proof. $\square$

PROOF OF PROPOSITION 3. When $M = 2N^{5/6}$, Proposition 5 gives

$$G_I^N(0,0) = 1 + \frac{M}{\sigma^2 N} + O(N^{-1/3} \log N) = 1 + O(N^{-1/6}),$$

and if $x = N^{5/6} + O(N^{1/2} \log N)$,

$$G_I^N(x,0) = -\frac{N^{5/6} + O(N^{1/2} \log N)}{\sigma^2 N} + \frac{2N^{5/6}}{\sigma^2 N} + O(N^{-1/3} \log N)$$

$$= \frac{1}{\sigma^2 N^{1/6}} + O(N^{-1/3} \log N).$$

Let $H_I^N(x,0)$ denote the probability that the random walk $X_k^N$ started at $x$ hits 0 before leaving $I$. Breaking things down according to the hitting time of 0:

$$G_I^N(x,0) = H_I^N(x,0) G_I^N(0,0)$$

which gives

$$H_I^N(x,0) = \frac{1}{\sigma^2 N^{1/6}} + O(N^{-1/3} \log N),$$

which is the desired result. □

LEMMA 5. *If* $0 < |x| < 2N^{5/6}$,

$$H_I^N(x,0) \leq \frac{2}{\sigma^2 N^{1/6}} + CN^{-1/3} \log N.$$

PROOF. Taking $y = 0$ in Proposition 5 we see that

$$G_I^N(x,0) \leq \frac{M}{\sigma^2 N} + CN^{-1/3} \log N.$$

The result now follows from $H_I^N(x,0) = G_I^N(x,0)/G_I^N(0,0)$. □

**7. Proof of Proposition 4.** The proof relies on a local central limit theorem, with bounds that take into account the dependence on $N$. First, we need a few definitions. Let

$$\rho_\ell(x) = \frac{1}{\sqrt{2\pi\sigma^2 \ell}} e^{-x^2/2\sigma^2 \ell}$$

be the normal density with variance $\ell\sigma^2$. Let $p_k^N$ be the distribution of the random walk at time $k$ when it starts at 0.

PROPOSITION 6 (Local central limit theorem). *Given a sequence of random walks with jump probabilities $p^N$ satisfying assumptions 1–5, there is a constant $C$, independent of $N$, such that for all $k \geq 1$ and all $x$ we have*

$$|p_k^N(x) - \rho_{kN}(x)| \leq \frac{C}{\sqrt{N} k^{3/2}}.$$

The proof of this uses standard techniques but is rather lengthy so we begin by giving the

PROOF OF PROPOSITION 4. By translation invariance it is enough to compute $a^N(x) = a^N(0, x)$. The local central limit theorem shows that, for all $N$ and $k \geq 1$,

$$p_k^N(0) - p_k^N(x) = \rho_{kN}(0) - \rho_{kN}(x) + O\left(\frac{1}{\sqrt{N} k^{3/2}}\right).$$

Therefore, after summing over $k \geq 0$,

$$a^N(x) = \sum_{k=0}^{\infty} p_k(x) - p_k(0)$$

$$= -1 + \delta(x, 0) - \sum_{k=1}^{\infty} [\rho_{kN}(0) - \rho_{kN}(x)] + O(1/\sqrt{N}).$$

We will now show that $\sum_{k=1}^{\infty} [\rho_{kN}(0) - \rho_{kN}(x)] = |x|/\sigma^2 N + O(1/\sqrt{N})$. Let $z = x/\sqrt{\sigma^2 N}$. Recalling the definition of $\rho_{kN}$,

$$\sum_{k=1}^{\infty} [\rho_{kN}(0) - \rho_{kN}(x)] = \frac{1}{\sqrt{2\pi\sigma^2 N}} \sum_{k=1}^{\infty} \frac{1}{\sqrt{k}} [1 - e^{-z^2/2k}].$$

Now the function $f_z(t) = (1 - e^{-z^2/2t})/\sqrt{2t}$, being a decreasing function divided by an increasing function, is decreasing in $t$ and therefore

$$0 \leq f_z(k) - \int_k^{k+1} f_z(t)\, dt \leq f_z(k) - f_z(k+1)$$

and thus $\sum_{k=1}^{\infty} f_z(k) - \int_1^{\infty} f_z(t)\, dt \leq f_z(1) \leq 1$. For the missing first piece of the integral we note

$$\int_0^1 f_z(t)\, dt = \int_0^1 \frac{1}{\sqrt{t}} (1 - e^{-z^2/2t})\, dt \leq \int_0^1 \frac{1}{\sqrt{t}}\, dt = 2.$$

Hence,

$$\sum_{k=1}^{\infty} [\rho_{kN}(0) - \rho_{kN}(x)]$$

$$= \frac{1}{\sqrt{\sigma^2 N}} \left( \int_0^{\infty} \frac{1}{\sqrt{2\pi t}} [1 - e^{-z^2/2t}]\, dt + O(1) \right).$$

We have put the $\sqrt{2\pi}$ inside so that the integral is $-1$ times the recurrent potential kernel for one-dimensional Brownian motion and hence is equal to $|z|$. One can find this fact on page 103 in Durrett [7], or derive it by changing variables $t = z^2/u$ and doing some calculus. In either case the result is

$$= \frac{1}{\sqrt{N}}\left(\frac{|z|}{\sigma} + O(1)\right) = \frac{|x|}{\sigma^2 N} + O\left(\frac{1}{\sqrt{N}}\right).$$

Putting everything together we get

$$a^N(x) = -1 + \delta(x,0) - \frac{|x|}{\sigma^2 N} + O\left(\frac{1}{\sqrt{N}}\right),$$

which is the desired result. □

Before entering into the proof of Proposition 6, we begin with an estimate on $\phi^N$, the characteristic function of the displacement $\xi_1^N$. This is the only proof that will require the use of assumption 3.

LEMMA 6.  *Let $L = N^{1/2}$. There are constants $a, b > 0$ so that, for all $N$ and $|\theta| \leq \pi$, we have*

$$|\phi^N(\theta)| \leq \begin{cases} (1 - bN\theta^2), & |\theta| \leq 4/(2L+1), \\ (1 - a), & |\theta| \in (4/(2L+1), \pi]. \end{cases}$$

*Consequently, given any $\varepsilon \in (0, \pi]$, there is a $c > 0$, independent of $N$, such that whenever $|\theta| > \varepsilon/\sqrt{N}$, $|\phi^N(\theta)| \leq e^{-c}$.*

PROOF. We use assumptions 3 to write $p^N(x) = b^N(x) + r^N(x)$ where $b^N(x) = h/L$ for $|x| \leq L$, 0 otherwise, and $r^N(x) = p^N(x) - b^N(x) \geq 0$. Define $B^N = \sum_x b^N(x) = (2L+1)h/L \to 2h$ as $L \to \infty$, and note that $\sum_x r^N(x) = 1 - B^N$. To bound $\phi^N$,

$$|\phi^N(\theta)| \leq \sum_{x=-\infty}^{\infty} r(x) + \left|\sum_{x=-L}^{L} b(x)e^{ix\theta}\right|$$

(17)
$$= (1 - B^N) + \frac{h}{L}\left|\frac{e^{i(L+1)\theta} - e^{-iL\theta}}{e^{i\theta} - 1}\right|$$

$$= (1 - B^N) + \frac{h}{L}\frac{|\sin((2L+1)\theta/2)|}{|\sin(\theta/2)|},$$

where in the last step we have multiplied numerator and denominator by $e^{-i\theta/2}$. Since the last expression is symmetric in $\theta$, we now restrict ourselves to $\theta \in [0, \pi]$.

For the next step we need the following inequalities.

LEMMA 7. (i) $\sin(\theta/2) > \theta/4$ for all $\theta \leq \pi$.
(ii) If $x > 2$, $(\sin 2)/2 > (\sin x)/x$.

PROOF. First we observe that $\cos x$ is decreasing on $[0, \pi/2)$ so if $x \leq \pi/2$, then

$$\frac{\sin x}{x} = \frac{1}{x}\int_0^x \cos y\, dy \geq \frac{2}{\pi}\int_0^{\pi/2} \cos y\, dy \geq \frac{1}{2},$$

which proves (i). For the second we note that

$$\left(\frac{\sin x}{x}\right)' = \frac{x\cos x - \sin x}{x^2}.$$

On $[2, \pi)$ the latter is negative since $\sin x > 0$ while $\cos x < 0$ there. Thus $\sin 2/2 > \sin x/x$ for all $x$ on this interval. For $x \in [\pi, 2\pi)$ the same inequality is obvious since $\sin x < 0$. Finally, for $x > 2\pi$, $\sin x/x < 1/x < 1/(2\pi) < 1/4 < \sin 2/2$, where we have used the fact that $\sin 2 \approx 0.909 > 1/2$. $\square$

Taking $x = (2L+1)\theta/2$, the inequalities in Lemma 7 imply that for $|\theta| > 4/(2L+1)$

$$\frac{h}{L}\frac{\sin((2L+1)\theta/2)}{\sin(\theta/2)} < \frac{h}{L}\frac{(\sin 2)(2L+1)\theta/4}{\theta/4} = B^N(\sin 2)$$

so we have

$$|\phi^N(\theta)| < 1 - (1 - \sin 2) \cdot B^N,$$

which gives the conclusion of Lemma 6 on this range of $\theta$.

For $\theta \leq 4/(2L+1)$, we rewrite the second term on the right-hand side of (17) as $Ee^{i\theta U}$ where $U$ is uniformly distributed on $\{-L, -L+1, \ldots, L\}$ and use (3.7) on page 101 of Durrett [8] to conclude

$$\left|Ee^{i\theta U} - \left(1 - \frac{\theta^2 EU^2}{2}\right)\right| \leq \frac{\theta^4}{4!}EU^4.$$

To bound the moments we use

$$EU^2 = \frac{2}{2L+1}\sum_{k=1}^L k^2 \geq \frac{2}{2L+1}\int_0^L x^2\, dx = \frac{2L^3}{3(2L+1)},$$

$$EU^4 \leq \int_{-L-1/2}^{L+1/2}\frac{x^4}{2L+1}\, dx = \frac{2(L+1/2)^5}{5(2L+1)}.$$

For $\theta \leq 4/(2L+1)$

$$\frac{\theta^2}{2}EU^2 \geq \theta^2\frac{N}{6}\cdot\frac{2L}{2L+1},$$

$$\frac{\theta^4}{4!}EU^4 \leq \frac{\theta^2}{4!} \cdot \frac{4^2}{(2L+1)^2} \cdot \frac{(L+1/2)^4}{5}$$

$$\leq \theta^2 \frac{N}{30}(1+1/2L)^2,$$

so we have

$$|Ee^{i\theta U}| \leq 1 - \frac{\theta^2 EU^2}{2} + \frac{\theta^4}{4!}EU^4$$

$$\leq 1 - \theta^2 \frac{N}{6} \cdot \frac{2L}{2L+1} + \theta^2 \frac{N}{30}(1+1/2L)^2.$$

Even when $L=1$, $(1/6) \cdot 2/3 = 1/9$ is larger than $(1/30) \cdot (3/2)^2 = 3/40$ and we have proved Lemma 6.  $\square$

PROOF OF PROPOSITION 6. By assumptions 2 and 4, we have that

(18) $$\phi^N(\theta) = 1 - \frac{\sigma^2 N}{2}\theta^2 + N^2 O(|\theta|^4).$$

The inversion formula gives

$$p_k^N(x) = \frac{1}{2\pi} \int_{-\pi}^{\pi} [\phi^N(\theta)]^k e^{-ix\theta} \, d\theta.$$

Introducing new variables $s = \sqrt{Nk} \cdot \theta$ and $z = x/\sqrt{Nk}$, the above is

$$= \frac{1}{2\pi\sqrt{Nk}} \int_{-\sqrt{Nk}\pi}^{\sqrt{Nk}\pi} \left[\phi^N\left(\frac{s}{\sqrt{Nk}}\right)\right]^k e^{-izs} \, ds.$$

Now, by (18), we can find an $\varepsilon > 0$, independent of $N$, such that if $s \leq \varepsilon\sqrt{Nk}$, the following approximation is valid:

$$\left[\phi^N\left(\frac{s}{\sqrt{Nk}}\right)\right]^k = \exp\left\{k \log \phi^N\left(\frac{s}{\sqrt{Nk}}\right)\right\}$$

$$= \exp\left\{k\left[-\frac{\sigma^2 s^2}{2k} + O\left(\left|\frac{s}{\sqrt{k}}\right|^4\right)\right]\right\}$$

$$= \exp\left\{-\frac{\sigma^2 s^2}{2}\right\} \exp(g(s,k)),$$

where $|g(s,k)| \leq c|s|^4/k$, and $c$ is independent of $N$. By choosing $\varepsilon$ smaller if necessary, we can guarantee that $|g(s,k)| \leq \sigma^2 s^2/4$.

As observed in Lemma 6, there is a $c$ independent of $N$ such that $|\phi^N(\theta)| \leq e^{-c}$ for $|\theta| > \varepsilon/\sqrt{N}$, that is, for $|s| > \varepsilon\sqrt{k}$. Using these observations we can

rewrite

$$p_k^N(x) = \frac{1}{2\pi\sqrt{Nk}} \int_{|s|\leq \varepsilon\sqrt{k}} e^{-izs} e^{-\sigma^2 s^2/2} e^{g(s,k)} \, ds$$
$$+ \frac{1}{2\pi\sqrt{Nk}} \int_{\varepsilon\sqrt{k}<|s|<\sqrt{Nk}\pi} \left[\phi^N\left(\frac{s}{\sqrt{Nk}}\right)\right]^k e^{-izs} \, ds.$$

The second term is $\leq e^{-ck}$ since the integrand is $\leq e^{-ck}$ by Lemma 6 and the interval has length $\leq 2\pi\sqrt{Nk}$.

Using $|e^{g(s,k)} - 1| < cs^4/k$ if $|s| \leq k^{1/4}$ and $|e^{g(s,k)} - 1| < \exp(\sigma^2 s^2/4)$ if $k^{1/4} < |s| \leq \varepsilon\sqrt{k}$, we have

$$\frac{1}{2\pi\sqrt{Nk}} \int_{|s|\leq \varepsilon\sqrt{k}} e^{-\sigma^2 s^2/2} |e^{g(s,k)} - 1| \, ds$$
$$\leq \frac{1}{2\pi\sqrt{Nk}} \left( \int_{-k^{1/4}}^{k^{1/4}} e^{-\sigma^2 s^2/2} cs^4/k \, ds + \int_{|s|\geq k^{1/4}} e^{-\sigma^2 s^2/4} \, ds \right).$$

Replacing $k^{1/4}$ by $\infty$ in the limits in the first integral, and using $s^2 \geq |s|k^{1/4}$ in the second, the above is

$$\leq \frac{1}{2\pi\sqrt{Nk}}\left(\frac{c}{k} + O(e^{-ck^{1/2}})\right).$$

On the other hand, setting $z = x/\sqrt{Nk}$ and later $s = \theta\sqrt{Nk}$ we have

$$\frac{1}{2\pi\sqrt{Nk}} \int_{|s|\leq \varepsilon\sqrt{k}} e^{-izs} e^{-\sigma^2 s^2/2} \, ds$$
$$= \frac{1}{2\pi\sqrt{Nk}} \int_{|s|\leq \varepsilon\sqrt{k}} e^{-ixs/\sqrt{Nk}} e^{-\sigma^2 s^2/2} \, ds$$
$$= \frac{1}{2\pi} \int_{-\infty}^{\infty} e^{-ix\theta} e^{-\sigma^2 Nk\theta^2/2} \, ds$$
$$- \frac{1}{2\pi\sqrt{Nk}} \int_{|s|>\varepsilon\sqrt{k}} e^{-izs} e^{-\sigma^2 s^2/2} \, ds$$
$$= \rho_{Nk}(x) + O(e^{-ck}),$$

which proves the result. □

## APPENDIX

Maruyama [11] considered a discrete-time Wright–Fisher model with a ring of $2n$ colonies with $N$ diploid individuals, nearest-neighbor migration with probability $m$ and mutation rate $u$ per generation. Here, to facilitate comparison with Maruyama [11] we use his notation. In Section 6 he considered sampling two individuals, one from the colony at 0 and the other at $i$, and let $f_i$ be the probability the two were identical by descent. This occurs if there is no mutation before the coalescence time $t_0$ so

$$f_i = E_i(1-2u)^{t_0}.$$

By writing recursive equations for the $f_i$ and then finding all of the eigenvalues and eigenvectors of an associated matrix, he developed exact but somewhat cumbersome formulas for the $f_i$. Writing $I_0$ instead of his $T$, (6.7) says

$$f_0 = \frac{(1-u)^2}{(1-u)^2 + 2N/I_0},$$

where $u$ is the mutation probability per generation,

$$I_0 = \pi^{-1} \int_0^\pi \frac{[1-m(1-\cos\theta)]^2}{1-(1-u)^2[1-m(1-\cos\theta)]^2} \, d\theta$$

and $m$ is the migration probability per generation. We are interested in the limiting behavior as $u \to 0$, so $I_0 \to \infty$. Since the contribution to the integral over $[\varepsilon, \pi]$ stays bounded, it is enough to investigate the behavior near 0. Since $1 - \cos\theta \sim \theta^2/2$ as $\theta \to 0$, the denominator should be well approximated by

$$1 - (1-2u)(1-m\theta^2) = 2u + m\theta^2.$$

Note that here $u$ and $\theta$ are small but $m$ need not be. Changing variables $\theta = (2u/m)^{1/2} x$, we have

$$\pi I_0 \sim \int_0^\infty \frac{1}{2u + 2ux^2} \left(\frac{2u}{m}\right)^{1/2} dx$$
$$= \frac{1}{(2um)^{1/2}} \int_0^\infty \frac{1}{1+x^2} \, dx$$
$$= \frac{\pi/2}{(2um)^{1/2}}$$

where to evaluate the integral we have used the definition of the Cauchy distribution (see, e.g., page 43 of Durrett [8]). Combining our calculations gives

$$f_0 \approx \frac{1}{1 + 4N(2um)^{1/2}}.$$

Changing notation $m = \nu$, and setting $u = \lambda/2(L^2/\nu)$, we have

$$E_0(\exp(-\lambda t_0/(L^2/\nu))) = f_0 \approx (1 + 4\sqrt{\lambda} N\nu/L)^{-1},$$

from which (1) follows.

**Acknowledgments.** The authors would like to thank Philip Protter and Dennis Talay for pointing out formula (2) in the *Handbook of Brownian Motion*, and thank Greg Lawler for useful discussions concerning the potential kernel of one-dimensional random walks and for the derivation of the boundary conditions of the PDE. The reports of two anonymous referees were useful in clarifying the exposition.

## REFERENCES

[1] ALDOUS, D. and FILL, J. A. (2002). Chapter 3 in *Reversible Markov Chains and Random Walks on Graphs*. Available at http://www.stat.berkeley.edu/users/aldous/RWG/book.html.

[2] BILLINGSLEY, P. (1971). *Weak Convergence of Measures*: *Applications in Probability*. SIAM, Philadelphia. MR0310933

[3] BORODIN, A. N. (1981). On the asymptotic behavior of local times of recurrent random walks with finite variance. *Theory Probab. Appl.* **26** 758–772. MR0636771

[4] BORODIN, A. N. and SALMINEN, P. (2002). *Handbook of Brownian Motion—Facts and Formulae.* Birkhäuser, Basel. MR1912205

[5] BOWEN, B. W. and GRANT, W. S. (1997). Phylogeography of the sardines (*Sardinops* spp.): Assessing biogeographic models and population histories in temperate upwelling regions. *Evolution* **51** 1601–1610.

[6] COX, J. T. and DURRETT, R. (2002). The stepping stone model: New formulas expose old myths. *Ann. Appl. Probab.* **12** 1348–1377. MR1936596

[7] DURRETT, R. (2005). *Stochastic Calculus*, 3rd ed. CRC Press, Boca Raton, FL. MR1398879

[8] DURRETT, R. (2004). *Probability*: *Theory and Examples*, 3rd Ed. Duxbury Press, Belmont, CA. MR1609153

[9] ITÔ, K. and MCKEAN, H. P. (1974). *Diffusion Processes and Their Sample Paths*, 2nd ed. Springer, Berlin. MR0345224

[10] JACOD, J. and SHIRYAEV, A. N. (1987). *Limit Theorems for Stochastic Processes*. Springer, Berlin. MR0959133

[11] MARUYAMA, T. (1970). Stepping stone models of finite length. *Adv. in Appl. Probab.* **2** 229–258. MR0266341

[12] MARUYAMA, T. (1971). The rate of decrease of heterozygosity in population occupying a circular or linear habitat. *Genetics* **67** 437–454. MR0452817

[13] MATSEN, F. A. and WAKELEY, J. (2005). Convergence to Island model coalescent process in populations with restricted migration. *Genetics* published online October 11, 2005.

[14] WILKINS, J. F. and WAKELEY, J. (2002). The coalescent in a continuous, finite, linear population. *Genetics* **161** 873–888.

[15] ZÄHLE, I., COX, T. and DURRETT, R. (2005). The stepping stone model. II. Genealogies and the infinite sites model. *Ann. Appl. Probab.* **15** 671–699. MR2114986

DEPARTMENT OF MATHEMATICS
CORNELL UNIVERSITY
523 MALOTT HALL
ITHACA, NEW YORK 14853
USA
E-MAIL: rtd1@cornell.edu

DEPARTMENT OF APPLIED MATHEMATICS
CORNELL UNIVERSITY
RHODES HALL
ITHACA, NEW YORK 14853
USA
E-MAIL: mr324@cornell.edu